\documentclass{amsart}

 \usepackage{amsfonts}

\begin{document}


\title[Generators of a toric ideal]
{Bounding The Degrees of Generators of a Homogeneous Dimension 2 
Toric Ideal}

\author{Hugh Thomas}

\address{Department of Mathematics\\University of Western Ontario
\\London, Ontario\\N6A 5B7 Canada}

\email{hthomas2@uwo.ca}

\thanks{Hugh Thomas is an Imperial Oil Post-doctoral 
Fellow at the University of Western Ontario}

\newcommand{\N}{\mathbb{N}}
\newcommand{\Z}{\mathbb{Z}}
\newcommand{\Tor}{\operatorname{Tor}^{R}}
\renewcommand{\Sigma}{\Lambda}
\newcommand{\Cal}[1]{\mathcal{#1}}
\newcommand\Sum{\sum}
\newtheorem{lemma}{Lemma}
\newtheorem{theorem}{Theorem}
\newtheorem{nonumalgorithm}{Algorithm}

\begin{abstract} Let $I$ be the toric ideal 
defined by a $2\times n$ matrix of integers, 
$$\Cal{A} = \left(\begin{array}{cccc}   1&   1 & \dots & 1 
\\ a_1 & a_2 & \dots & a_n
\end{array}\right)$$ 
with $a_1 < a_2<\dots <a_n$.  We give a combinatorial proof that $I$
is generated by elements of degree 
at most the sum of the two largest differences $a_{i}-a_{i-1}$.  
The novelty is in the method of proof: the result has already been shown by 
L'vovsky using cohomological arguments. \end{abstract}



\maketitle

\section*{Introduction}

Let $\Cal{A}=(a_{ij})$ be a $d \times n$ matrix of integers.  Let $k$ be an arbitrary
ground field.  Let $R=k[z_1,\dots,z_n]$ and $S=k[y_1,\dots,y_d]$.  From 
$\Cal{A}$ we get a ring homomorphism
 $\psi$ from $R$ to $S$ by sending $z_j$ to 
$y_1^{a_{1j}}y_2^{a_{2j}}\dots y_d^{a_{dj}}$.  Let $I$ be the kernel of this map.  Ideals which arise in this way are called toric ideals.  
See \cite{St} for a thorough introduction
to the subject.  

It is natural to try to determine the syzygies of such an ideal $I$,
a problem pursued in \cite{CM,CP,CG2}, 
or, more restrictedly, to ask for a minimal set of generators
for such an ideal, the approach taken in \cite{BCMP}.  
In our case, as in \cite{BCP}, 
we shall be interested in determining an upper bound for the degrees of
a minimal generating set, in a special case, 
also singled out for consideration in \cite{CG, Lv, BCP}, 
as follows:
let $d=2$, and let
all the $a_{1j}=1$.  (It follows that the ideal $I$ is the homogeneous 
ideal of a monomial curve
in projective space, but we shall not adopt that point of view here.)
For simplicity, we refer to $a_{2j}$ as 
$a_j$.  Without loss of generality let  $a_1<a_2<\dots<a_n$.  

The ring $R$ has a $\N \times \Z$ grading, where $z_j$ has degree $(1,a_{j})$.
It is easily seen that the ideal $I$ is homogeneous with respect to this 
grading.  Forgetting the $\Z$ component of the grading, we recover the 
usual $\N$ grading on $R$.  We use the word ``bidegree'' to refer to 
degree in the $\N \times \Z$ grading, and ``degree'' to refer to degree in 
the usual $\N$ grading.

This paper consists of a proof of the following theorem:

\begin{theorem}[Main Theorem]
Let $r\geq s$ be the two largest successive
differences $a_i-a_{i-1}$. 
Then $I$ is generated by elements of degree no more than $r+s$.

\end{theorem}

Simple examples show this bound is tight: let $r$ and $s$ be relatively prime
integers, let
$n=3$, $a_1=-r$, $a_2=0$, $a_3=s$.  Then $I$ is the principal ideal generated
by $z_2^{r+s}-z_1^sz_3^r$.  

Note that a given toric ideal $I$ will arise from more than one choice of
$\Cal{A}$.  Given a toric ideal $I$ in $R$, 
to use the Main Theorem to obtain the best possible bound for the
degrees of a minimal generating set for $I$, one should choose $\Cal{A}$ in
such a way that the greatest common divisor of the successive differences
$a_i-a_{i-1}$ is 1.  Such a choice is always possible, and any such choice
will yield the same bound.  

In \cite{Lv}, L'vovsky used cohomological results from \cite{GLP}
to prove a stronger statement than our Main Theorem,
bounding the regularity of $I$, which is to say, bounding the degrees of
the $i$-syzygies of $I$ for all $i$.  Restricted to generators (0-syzygies),
his result coincides with ours.  It is not clear how the bound given
by our Main Theorem compares with the bound obtained in \cite{BCP}.

In justification for this paper, 
aside from the intrinsic interest of a combinatorial proof of L'vovsky's 
bound,
we hope
that the techniques
of this proof may extend to higher dimensional cases.  

The Main Theorem follows easily from the following combinatorial result.

\begin{theorem}[Connectedness Theorem] Let $V \subset \Z$, not necessarily finite, 
with the sizes of gaps between successive elements bounded above.
Let $r \geq s$ be the sizes of the two largest gaps 
between successive elements of
$V$.  
For $(q,c)\in \N \times \Z$, let $\Pi_{(q,c)}$ be the collection of 
multisets with support in $V$, of cardinality $q$ and sum $c$.  Let 
$\Delta_{(q,c)}$ be the simplicial complex generated by the supports of the
multisets in $\Pi_{(q,c)}$.  Then if $q>r+s$, $\Delta_{(q,c)}$ is connected.
\end{theorem}

\section*{Translation to Combinatorics}

A multiset is an unordered collection of elements, in which some elements
may appear with multiplicity greater than one.
We use $+$ and $-$ for addition and subtraction of multisets, and 
write $P=\{x_1,\dots,x_n\}_\leq$ to indicate that the elements of $P$ are
listed in non-decreasing order.  $\sum P$ is the sum of the elements of $P$.
We now give the proof of the Main Theorem assuming the Connectedness 
Theorem:
\begin{proof}[Proof of Main Theorem] Let $V=\{a_1,\dots,a_n\}_<$.  
Let $r\geq s$ be the
two largest successive differences $a_i-a_{i-1}$.  By the Connectedness
Theorem, $\Delta_{(q,c)}$ is connected for $q>r+s$.  The Main Theorem
now follows from the following lemma:

\begin{lemma}[Translation Lemma] No minimal 
generating set for $I$ has generators
in bidegree $(q,c)$ iff $\Delta_{(q,c)}$ is connected.  
\end{lemma}

Remark: 
This is a special case of a result of \cite{CP}, which gives
information about the degrees of minimal generators of $I$ and also 
of all its $i$-syzygies, based on the homology of $\Delta_{(q,c)}$.
In the interest
of self-containedness, we give an elementary proof of the result we need.

\begin{proof}[Proof] ($\Leftarrow$)   
Let $I^<$ denote the ideal of $R$ generated by 
the elements of $I$ of bidegree $(q',c')$ with $q'<q$.  
Then $I^<_{(q,c)}$, the $(q,c)$-bigraded part of $I^<$, is a sub-vector space
of $I_{(q,c)}$.  
We wish to show that
$I_{(q,c)}=I^<_{(q,c)}$.

For $P$ a multiset with support in $V$, let $z^P$ denote
the monomial in $R$ where the exponent of $z_i$ is the multiplicity of 
$a_i$ in $P$.  
$I_{(q,c)}$ is spanned as a $k$-vector space by elements of the form
$z^P-z^{P'}$, for $P$ and $P'$ in $\Pi_{(q,c)}$.  

Let $P$ and $P'$ be two elements of $\Pi_{(q,c)}$ with non-empty intersection,
say $A$.  Then $z^P-z^{P'}=z^A(z^{P- A}-z^{P' - A})$.
But $z^{P - A}-z^{P'- A}$ is in $I$, in bidegree
$(|P- A|,\sum (P-A))$.  Thus $z^P-z^{P'} \in I^<_{(q,c)}$.  

Now suppose $P$ and $P'$ are arbitrary elements of $\Pi_{(q,c)}$.  
Since $\Delta_{(q,c)}$ is connected, it follows that there exist $P=P_0,
P_1,\dots,P_t=P'$ with each $P_i \in \Pi_{(q,c)}$, such that for each
$i$, $P_i$ and $P_{i+1}$ have at least one 
element in common.  By the previous argument,
then, $z^{P_{i}}-z^{P_{i+1}} \in I^<_{(q,c)}$, from which it follows that
$z^P-z^{P'} \in I^<_{(q,c)}$, as desired.

($\Rightarrow$) $I^<_{(q,c)}$ is spanned by elements of the form
$z^A(z^B-z^{B'})$ with $A$ non-empty, or equivalently of the form $z^P-z^{P'}$ 
with $P$ and $P'$ having non-empty intersection.  
 It follows that for any $z^P-z^{P'}
\in I^<_{(q,c)}$, the supports of $P$ and $P'$ are in the same component of
$\Delta_{(q,c)}$.  Thus, if $\Delta_{(q,c)}$ has more than one component,
$I^<_{(q,c)} \ne I_{(q,c)}$, as desired.  
\end{proof} \end{proof}

\section*{Combinatorial Lemmas}

We now develop the combinatorial tools to prove the Connectedness Theorem.
Fix $V \subset \Z$, with $r$ and $s$ the sizes of the two largest gaps between 
successive elements of $V$.  
Fix $q>r+s$ and fix $c\in \Z$.  
Let $\Delta=\Delta_{(q,c)}$ and $\Pi=\Pi_{(q,c)}$.  

Given
a multiset                         $P=\{x_1,\dots,x_p\}_\leq$, define $m(P)=
\sum_i ix_i$.  Among multisets of the same sum and cardinality, 
the intuition is 
that $m$ is a measure of how
spread out $P$ is --- the more spread out $P$ is, the greater $m(P)$ will be.  

We now show that, for a suitable class of $P \in \Pi$, we can find
another multiset $P' \in \Pi$ which is more spread out than $P$.  

\begin{lemma}[Expansion Lemma] 
Let $P \in \Pi$, $P=A+C$.  $|C|=r+s$, and $C$ contains neither the greatest 
element nor the least element of $V$.  Then there exists some 
$P' \in \Pi$, $P'=A+C'$, such that $m(P')>m(P)$.
\end{lemma}

Note: the following proof owes its basic approach 
to the proof of Theorem 6.1 of \cite{St}.

\begin{proof}[Proof]
Consider the following algorithm, which obtains a sequence of multisets 
$C_i$, with $C_0=C$, where $C_i$ is obtained from $C_{i-1}$ by replacing one
of the original elements of $C$ by either the next larger or the next 
smaller element of $V$.  When thinking about this algorithm, it's helpful
to think of the elements of $C$ as stones sitting on a number line, where
the allowed positions for the stones are the numbers in $V$.  $C_i$ is 
obtained from $C_{i-1}$ by jumping one stone which hasn't been moved yet 
to the next higher or lower
allowed position.  

\begin{nonumalgorithm}[Expansion Algorithm]\end{nonumalgorithm}
\noindent
$C_0 := C$. \\ 
$Active := C$.  (These are element of $C$ which haven't moved yet.) \\
$s_0 := 0$. \\
For $i :=1$ to $r+s$ do: \\ \hspace*{.75cm}
    If $s_{i-1} \leq 0$, \\ \hspace*{1.5cm}
        Remove the largest element from Active, and call it $x_i$. \\ 
\hspace*{1.5cm}
        $C_{i}:=C_{i-1}-\{x_i\}+\{\mbox{the least element of $V$ greater than
 $x_i$}\}$. \\ \hspace*{.75cm}
    If $s_{i-1} > 0$, \\ \hspace*{1.5cm}
        Remove the smallest element from Active, and call it $x_i$. 
\\ \hspace*{1.5cm}
        $C_{i}:=C_{i-1}-\{x_i\}+\{\mbox{the greatest element of $V$ less
 than $x_i$}\}$.  \\ \hspace*{0.75cm}
   Let $s_i := \Sum C_i-\sum C$.  

\bigskip

We would now like to bound $s_i$.  Observe that if $s_{i-1}\leq 0$ then
$s_i>s_{i-1}$, and if $s_{i-1}>0$, then $s_i<s_{i-1}$.
Thus, the absolute value of $s_i$ can be no greater than the largest jump
possible on a single step, which is $r$.  However, we can be a little more 
precise.  

The elements of $C$ which are increased by the algorithm are greater than
or equal to all the elements which are decreased by the algorithm.  Thus,
any gap between successive elements of $V$ can be jumped in only one direction
in the course of running this algorithm (though it may be jumped more than
one time).  

Suppose there is a unique gap of size $r$, and it is jumped only
in the increasing
direction.    Then it follows that $-s+1 \leq s_i \leq r$ for all $i$.  
Symmetrically, suppose there
is a unique gap of size $r$, and it is jumped only in the decreasing 
direction.  Then $-r+1 \leq s_i \leq s$.  If there is a unique gap of size
$r$ which is not jumped, or there is more than one gap of size $r$ (in which 
case $r=s$), then $-s+1 \leq s_i \leq s$.  In any case, we deduce that there
are at most $r+s$ possible values for $s_i$.  But there are $r+s+1$ of the
$s_i$.  Thus, at least two of the $s_i$ must be equal, say $s_j=s_l$, with
$j<l$.  Obtain $C'$ from $C$ by making the same jumps as were made
in the algorithm on steps $j+1$ through $l$.  Then $\Sum C' = \Sum C +s_l-
s_j=\Sum C$.  Let $P'=A+C'$.  In going from $P$ to $P'$, the elements
which have been increased are all greater than or equal to the elements which
have been decreased, and thus $m(P')>m(P)$.  
\end{proof}

Using only this lemma, we can prove the Connectedness Theorem with the 
additional assumption that $V$ is not bounded below.

\begin{proof}[Proof of Connectedness Theorem assuming $V$ is not bounded below]  
Let $x$ and $y$ be vertices
of $\Delta$.  We want to show that they lie in the same component.
Choose some $P \in \Pi$ with  $x \in P$, and some 
$P' \in \Pi$, with $y\in P'$.  Let $x'$ be the maximum element of $P$,
$y'$ the maximum element of $P'$.  If $x'=y'$ then we are done. 
So assume without loss of generality 
that $x'<y'$.  Choose a subset $A$ of $P$, whose size is 
$q-(r+s)$, and which contains $x$.  Thus, we can write $P=A+C$, satisfying
the hypotheses of the Expansion Lemma.  

Now apply the Expansion Lemma recursively, with this choice of $A$ fixed,
obtaining a sequence of multisets $P=P_0,P_1,P_2,\dots,P_t$ until either
$m(P_t) \geq (1+2+\dots+q)y'$ or $P_t$ contains the largest element of $V$.  
In either case, $P_t$ clearly contains
an element which is greater than or equal to $y'$.  Thus, some previous 
$P_i$ contains $y'$.  Now by construction $P_i \in \Pi$, and $x$ and $y'$
are both contained in $P_i$, which finishes the proof.  
\end{proof}

\begin{lemma}[Multiple Expansion Lemma] Let $P, P' \in \Pi$.  Suppose the largest
and smallest elements of $P+P'$ occur in $P'$.  Then there is a $Q\in \Pi$
containing at least one element of each of $P$ and $P'$.  \end{lemma}

\begin{proof}[Proof] The proof is essentially the argument given above, 
proving the
Connectedness Theorem in the case where $V$ is not bounded below (but we do
not make the assumption that $V$ is not bounded below).
Fix some set 
$A$ in $P$ of size $q-(r+s)$.  Apply the Expansion Lemma recursively with
$A$ fixed, producing a sequence of multisets $P=P_0, P_1,\dots$, until 
some $P_t$ includes some element of $P'$.  The Expansion Lemma can always be
applied because at no stage before halting does $P_i$ include the largest 
or smallest element of $V$, since prior to including one of these elements,
it would include an element of $P'$.  Then $Q=P_t$ satisfies the conditions
in the statement of the lemma.  
\end{proof}

We now prove another lemma, similar to the Multiple Expansion Lemma, which we
will need to prove the Connectedness Theorem in full generality.  

\begin{lemma}[Criss-Cross Lemma] Let $P, P' \in \Pi$.  
Suppose the largest element of $P+P'$ occurs in $P'$, 
while the smallest element of $P+P'$ occurs in $P$.  
Then there exists a $Q \in \Pi$
which contains at least one element of each of $P$ and $P'$.  
\end{lemma}

\begin{proof}[Proof]
Assume, without loss of generality, that $P$ and $P'$ are disjoint.   
Pick $f \in P$, $f' \in P'$, with $f<f'$.  Let $B = P - \{ f\}$, 
$B'=P' - \{ f'\}$.

Split $B$ into two multisets, $X$ and $Y$, where $X$ consists of 
the elements of $B$ less than all elements of $B'$, and $Y$ is the remainder.
Similarly, split $B'$ into $X'$, the elements greater than all elements of 
$B$, 
and $Y'$, the remainder.  

\begin{lemma}[Size Lemma] Either $|Y|$ is greater than the longest gap below $B'$
or $|Y'|$ is greater than the longest gap above $B$.  \end{lemma}

\begin{proof}[Proof]
First, 
I claim that $|Y|+|Y'|>r+s$.  Let $B=\{b_1,\dots,b_{q-1}\}_\leq$, 
$B'=\{b'_1,\dots,
b'_{q-1}\}_\leq$.  For $1 \leq i \leq |X|$, and for 
$q-|X'| \leq i \leq q-1$, $b_i<b'_i$.  So, if $|X|+|X'| \geq q-1$, then
$\sum B < \sum B'$.  But $\sum B=\sum P -f > \sum P' -f' = \sum B'$, which 
is a contradiction.  Thus, $|X|+|X'|<q-1$, so
$|Y|+|Y'|>q-1\geq r+s$.  

Now, we can see that at least one of the following four cases holds:

\begin{enumerate}
\item There is a unique gap of size $r$ below $B'$ and $|Y| \geq r+1$.
\item All the gaps below $B'$ are of size no more than $s$ and $|Y| \geq s+1$.
\item There is a unique gap of size $r$ above $B$ and $|Y'| \geq r+1$.
\item All the gaps above $B$ are of size no more than $s$ and $|Y'|\geq s+1$.
\end{enumerate}

If there is a unique gap of size $r$ below $B'$ then either $|Y| \geq r+1$
(case (1))
or $|Y'| \geq s+1$ (case (4)).
Similarly, if there is a unique gap of size $r$ above $B$, we are in case
 (2) or case (3).
Otherwise, we are in case (2) or (4).  

It is now easy to see that in cases (1) or (2), $|Y|$ is greater than the
largest gap below $B'$, while in cases (3) or (4), $|Y'|$ is greater than
the largest gap above $B$, proving the lemma.   
\end{proof}

Now, using the Size Lemma, by symmetry, 
we may assume without loss of generality that
$|Y|$ is greater than the largest gap below $B'$.  Let the size of this gap
be $g$.    
Now consider the following algorithm.  

\begin{nonumalgorithm}[Criss-Cross Algorithm]\end{nonumalgorithm}

\noindent
$B_0:=B$.\\
$ActiveX:=X$. (These are the elements of $X$ that haven't moved yet.) \\
$ActiveY:=Y$. (And similarly for $Y$.) \\
$i:=0$. \\
$s_0:=0$. \\
Repeat until $s_i = s_j$ for some $j<i$: \\ \hspace*{0.75cm}
     $i:=i+1$. \\ \hspace*{0.75cm}
     If $s_{i-1} \leq 0$, then \\ \hspace*{1.5cm}
         If ActiveX is non-empty, \\ \hspace*{2.25cm}
             Remove the largest element of ActiveX and call it $x_i$.
\\ \hspace*{1.5cm}
         Otherwise, \\ \hspace*{2.25cm}
             Remove the largest element of ActiveY and call it $x_i$.
\\ \hspace*{1.5cm}
         $B_{i}:=B_{i-1}-\{x_i\}+\{\mbox{the least element of $V$ greater
 than $x_i$}\}$. \\ \hspace*{0.75cm}
     If $s_{i-1}>0$, \\ \hspace*{1.5cm}
         Remove the smallest element of ActiveY and call it $x_i$.
\\ \hspace*{1.5cm}
         $B_{i}:=B_{i-1}-\{x_i\}+\{\mbox{the greatest element of $V$ less
 than $x_i$}\}$. \\ \hspace*{0.75cm}
     Let $s_i=\sum B_i-\sum B$. 
\bigskip

I first argue that there are always enough elements of ActiveX and 
ActiveY to run this algorithm.  There are two potential problems:

First, suppose $s_{i-1}\leq 0$, $\text{ActiveX}=\emptyset$, $\text{ActiveY}=
\emptyset$.  Then $i \geq q$.  So at least $r+s+1$ values of $s_j$ have
been defined.  But, as in the proof of the Expansion Lemma, no gap can get 
jumped in both directions, so there are at most $r+s$ possible values for
the $s_j$, so the algorithm could never have reached this step, since 
there must 
already have been some value $s_j$ which has occured more than once.  

Second, suppose $s_{i-1}>0$, and $\text{ActiveY}=\emptyset$.  If 
$\text{ActiveX}=\emptyset$, then $i \geq q$, and we have a contradiction as in 
the previous case.  So assume that $\text{ActiveX}\ne \emptyset$.  

Thus, every element of $Y$ that have been removed from $\text{ActiveY}$ 
was removed on a turn $j$ with $s_j>0$.  Further, we can say that 
$g\geq s_j$ for all these $j$, since we have not yet jumped any elements of
$Y$ in the positive direction, so all the gaps that we  have jumped in the 
positive direction are no longer than $g$.  Since the algorithm has not
yet terminated, all the $s_j$ are different, and they have only $g$ 
different possible positive values.  Thus, the case $s_j>0$ has
been encountered at most $g$ times.  But since $|Y| \geq g+1$, it follows
that $\text{ActiveY} \ne \emptyset$, which contradicts the assumption
that we had run out of ActiveY.

Thus, the algorithm runs successfully, and we obtain $i>j$ with 
$s_i=s_j$.  Obtain $\tilde B$ from $B$, by making the same jumps as were made
on steps $j+1$ through $i$ of the algorithm.  
As in the proof of the Expansion Lemma, $\sum \tilde B = \sum B$.   

Let $P_0:=P$, $P_1:=\{f\}+\tilde B$.  If $P_1$ contains an element of 
$B'$, we are done with $Q=P_1$.  Otherwise, we apply the algorithm 
recursively, keeping  $B'$, $f$, and $f'$ fixed.  I claim 
that we eventually obtain a   $P_i$ which contains an element of $B'$.  

Note first that, until we obtain a $P_i$ which contains an element of $B'$,
we stay in the same case of the Size Lemma,
because no
element of $B$ can jump past an element of $B'$ without hitting it, 
so $X$ and $Y$
will stay the same size.  
Thus, we can succesfully repeat the algorithm until some $P_i$ 
intersects $B'$.

For a multiset $D=\{x_1\dots,x_p\}_\leq$, define $\tilde m(D) = 
\sum_i (p+1-i)x_i$.
This function plays a similar role to $m$, 
except that as a set gets more spread out, $\tilde m$
decreases.  

Now, I claim that $\tilde m( Y)$ is decreased each time we run the algorithm,
since the decreases to $Y$ are at least as great as its increases, and
the elements that are decreased precede the elements that are increased.

Thus, at some point, if the algorithm could run forever,
 $\tilde m(Y)$ would be smaller than $(1+2+\dots+|Y|)$ times the smallest
element of $B'$, but this would mean that $Y$ contained elements smaller than
every element of $B'$, which would be a contradiction.

So eventually, some $P_i$ contains an element of $B'$.  Then $Q=P_i$ satisfies
the conditions of the statement of the lemma.    
\end{proof}

\section*{Proof of the Connectedness Theorem}

Now, we prove the connectedness of $\Delta$ in full generality.

\begin{proof}[Proof of Connectedness Theorem] We begin as in the
case of $V$ not bounded below.  
Let $x$ and $y$ be vertices of $\Delta$.  
We want to show that $x$ and $y$ are in the same 
component.  Choose a set $P$ in $\Pi$ which contains $x$, and a set 
$P'$ in $\Pi$ which contains $y$.  

If $P$ and $P'$ intersect, then we are done.  So assume they do
not.  
Suppose that the greatest and the least elements of $P+P'$ both occur in 
one of $P$ or $P'$, without loss of generality, say $P'$.
Now apply the Multiple Expansion Lemma to obtain a multiset $Q \in \Pi$ 
intersecting
both $P$ and $P'$.  As in the proof of the Connectedness Theorem with $V$
not bounded below, this shows that $x$ and $y$ are connected in 
$\Delta$.  

Now, suppose that the greatest and least elements of $P+P'$ do not both 
occur in either $P$ or $P'$.  
Without loss of generality, 
let $P'$ contain the greatest.  This
puts us in the position to apply the Criss-Cross Lemma, and again, the 
multiset
$Q$ 
 which we obtain from it shows that $x$ and $y$ are connected in $\Delta$.  
\end{proof}

\section*{Further Directions}

First, it would be good to give a combinatorial proof of the entire result
of L'vovsky (bounding the degrees of all $i$-syzygies of $I$, not 
just the generators).  An argument might use the full strength of 
the result mentioned in the proof of the Translation Lemma from 
 \cite{CP} to translate the problem into a 
combinatorial framework. The necessary 
combinatorial result might then be established by an induction argument with
a version of the Connectedness Theorem as a base case, but so far we have
been unable to accomplish this 
except in the case where $V$ is not bounded below,
which is of limited interest.

Also, 
as mentioned in the introduction, one might apply the techniques
of this paper
to prove degree bounds or regularity bounds in higher dimensions.
However, this is considerably trickier.  Our strategy for modifying a 
multiset $P$ until it hits an element of $P'$ is akin to 
a game of hide-and-seek in one dimension; the interested reader will see
why the game is usually played in two dimensions.

\section*{Acknowledgments} 
I would like to thank Hal Schenck for suggesting this problem to me, and for 
fruitful initial discussions.  I would also like to thank Philippe Gim\'enez
and the referee
for their helpful comments.

\end{document}